\documentclass[11pt]{amsart}

\usepackage[latin1]{inputenc}
\usepackage[T1]{fontenc}
\usepackage{amsmath,amsfonts,amssymb,amsthm,url,geometry,mathrsfs,amsthm,graphicx}

\newcommand{\e}{\varepsilon}
\renewcommand{\phi}{\varphi}
\newcommand{\iy}{\infty}

\DeclareMathOperator{\diag}{diag} 
 
\DeclareMathOperator{\tr}{Tr}

\DeclareMathOperator{\Var}{\mathbf{Var}}
\DeclareMathOperator{\Cov}{\mathbf{Cov}}

\renewcommand{\leq}{\leqslant}
\renewcommand{\geq}{\geqslant}

\newcommand{\st}{\  : \ }
\newcommand{\N}{\mathbf{N}}

\newcommand{\R}{\mathbf{R}}
\newcommand{\C}{\mathbf{C}}
\DeclareMathOperator{\E}{\mathbf{E}}
\DeclareMathOperator{\card}{card}

\renewcommand{\P}{\mathbf{P}}
\newcommand{\M}{\mathcal{M}}

\newcommand{\mD}{\mathcal{D}}

\newcommand{\mM}{\mathcal{M}}

\newcommand{\Id}{\mathrm{Id}}

\newcommand{\ketbra}[2]{| #1 \rangle \langle #2 |}

\newcommand{\ket}[1]{| #1 \rangle}

\theoremstyle{plain}
\newtheorem*{thm}{Theorem}

\newtheorem*{conj}{Conjecture}

\newtheorem*{remark}{Remark}

\newtheorem{N-thm}{Theorem}
\newtheorem{N-cor}{Corollary}
\newtheorem{N-prop}{Proposition}[section]
\newtheorem{N-lem}[N-prop]{Lemma}
\newtheorem{N-fact}[N-prop]{Fact}

\theoremstyle{definition}
\newtheorem{N-defn}[N-prop]{Definition}
\newtheorem*{defn}{Definition}

\theoremstyle{remark}

\newtheorem*{eg}{Example}

\newcommand{\bra}[1]{\langle #1 |}

\title{Partial transposition of random states and non-centered semicircular distributions}
\author{Guillaume AUBRUN}
\address{Institut Camille Jordan, Universit\'e Claude Bernard Lyon 1, 43 boulevard du 11 novembre 1918, 69622 Villeurbanne CEDEX, France}
\email{aubrun@math.univ-lyon1.fr}

\thanks{This research was supported by the {\itshape Agence Nationale de la Recherche} grants ANR-08-BLAN-0311-03 and ANR-2011-BS01-008-01, and by the {\itshape Institut Mittag-Leffler} in Stockholm. I also thank S.~Szarek and I.~Nechita for useful comments.}

\begin{document}

\begin{abstract}
Let $W$ be a Wishart random matrix of size $d^2 \times d^2$, considered as a block matrix with $d \times d$ blocks. Let $Y$ be the matrix obtained by transposing each block of $W$. We prove that the empirical eigenvalue distribution of $Y$  approaches a non-centered semicircular distribution when $d \to \iy$. We also show the convergence of extreme eigenvalues towards the edge of the expected spectrum. The proofs are based on the moments method.

This matrix model is relevant to Quantum Information Theory and corresponds to the partial transposition of a random induced state. A natural question is: ``When does a random state have a positive partial transpose (PPT)?''. We answer this question and exhibit a strong threshold when the parameter from the Wishart distribution equals $4$. When $d$ gets large, a random state on $\C^d \otimes \C^d$ obtained after partial tracing a random pure state over some ancilla of dimension $\alpha d^2$ is typically PPT when $\alpha>4$ and typically non-PPT when $\alpha<4$. 
\end{abstract}

\maketitle

\section{Introduction} \label{sec:intro}

In the recent years, several connections were established between Random Matrix Theory and Quantum Information Theory. 
It turns out that random operators, and the random constructions they induce, can be used to construct quantum channels with an unexpected behavior, violating some natural conjectures (the most prominent example being Hastings's counterexample to additivity conjectures \cite{hastings}). Random matrices appear to be a sharp tool in order to understand the high-dimensional objects from Quantum Information Theory.

In this spirit, we study here a model of random matrices motivated by Quantum Information Theory. The model is simple to describe: start from Wishart $n \times n$ random matrices, which is the most natural model of random positive matrices. Assume that their dimension is a square ($n=d^2$). These matrices can be considered as block-matrices, with $d^2$ blocks, each block being a $d \times d$ matrix. Now our model is obtained by applying the transposition operation inside each block. A equivalent formulation is to consider $d^2 \times d^2$ matrices as operators on the tensor product of two $d$-dimensional spaces, and to apply to them the {\itshape partial transposition} $\Id \otimes T$, where $T$ is the usual transposition.

For this model, the empirical eigenvalue distribution converges towards a {\itshape non-centered} semicircular distribution, and the extreme eigenvalues converge towards the edge of the spectrum. These results were observed numerically by \v{Z}nidari\v{c} et al. \cite{znidaric}. The aim of the present paper is to give a complete proof of these facts. We rely on a standard tool from Random Matrix Theory: the method of moments.

The fact that the limiting distribution is semicircular is not a complete surprise. In the context of free 
probability, semicircular distributions are the non-commutative analogue of Gaussian distributions, and therefore one expects their appearance in limit theorems. For example Wigner's celebrated theorem identifies the centered semicircular distribution as the limit distribution of eigenvalues of random Hermitian matrices. However, other limiting distributions do appear in the theory: for example, the Wishart matrices themselves (i.e., {\em without} the partial transposition) converge to the so-called Mar\v{c}enko--Pastur law (see section
\ref{subsection:MP}). Moreover, our model brings some additional exoticism since the
limiting distribution is non-centered.

Since the transposition is not a completely positive map, there is no reason {\itshape a priori} for matrices 
from our model to be positive. However, we show that for some range of the parameters, partially transposed Wishart matrices are typically positive. A threshold occurs when the parameter from the Wishart distribution equals $4$.

The partial transposition appears to play a central role in Quantum Information Theory and is closely related to the concept of entanglement. An important class of states is the family of states with a Positive Partial Transpose (PPT). Non-PPT states are necessarily entangled \cite{peres} and this is the simplest test to detect entanglement. Let us simply mention a related important open problem known as the distillability conjecture \cite{hhh}: it asks whether, for a state $\rho$, non-PPT is equivalent to the existence of a protocol which, given many copies of $\rho$, distills them to obtain Bell singlets---the most useful form of entanglement. A positive answer to the distillability conjecture would give a physical meaning to partial transposition.

The model of Wishart random matrices has also a physical interpretation in terms of open systems: assume the subsystem $\C^d \otimes \C^d$ is coupled with some environment $\C^p$. If the overall system is in a random pure state, the state on $\C^d \otimes \C^d$ obtained by partial tracing over $\C^p$ is distributed as a (normalized) Wishart matrix. Early notable works about entanglement of random states include \cite{kendon} and \cite{hlw}. Our results can be translated in this language. In particular, a random induced state is typically non-PPT  when $p/d^2<4$ and is typically PPT when $p/d^2>4$. This shows that a threshold for the PPT property occurs at $p = 4d^2$.

\subsection*{Organization} The paper is organized as follows: Sections \ref{sec:background}--\ref{sec:extreme} are written in the language of Random Matrix Theory and contain the proof of our theorems. Section \ref{sec:background} introduces the model and states Theorem \ref{main-theorem} (convergence towards the non-centered semicircle distribution) and Theorem \ref{thm:extreme-egv} (convergence of the extreme eigenvalues). Section \ref{sec:wishart} reminds the reader about non-crossing partitions and the combinatorics behind the moments method for Wishart matrices, on which we rely heavily. Section \ref{sec:overview} shows how to derive Theorem \ref{main-theorem} from moment estimates ; the proof of these estimates (the heart of the moments method) is deferred to Sections \ref{sec:expectation} and \ref{sec:variance}. Section \ref{sec:extreme} contains the proof of Theorem \ref{thm:extreme-egv}. Section \ref{sec:QIT} connects to Quantum Information Theory. Section \ref{sec:misc} contains some general remarks and possible variations on the model. A high-level non-technical overview of the result of this paper and of a related article \cite{asy} can be found in \cite{asy-short}.

\section{Background and statement of the main theorem} \label{sec:background}

\subsection{Conventions} By the letters $C,C_0,c,\dots$ we denote absolute constants, whose value may change from occurrence to occurrence. The integer part of a real number $x$ is denoted by $\lfloor x \rfloor$. We denote by $[k]$ the set $\{1,\dots,k\}$. Addition in $[k]$ is understood modulo $k$. We denote by $\vec a,\vec b,\vec c,\dots$ multi-indices which
are elements of $\N^k$ for some integer $k$. The coordinates of $\vec a$ are denoted $(a_1,\dots,a_k)$. 

When $\vec a \in \N^k$, we denote by $\# \vec a$ the number of distinct elements which appear in the set $\{a_1,\dots,a_k\}$.  For example, $\#(1,4,1,2)=3$. The cardinality of a set $A$ is denoted $\card A$. The notation ${\bf 1}_E$ denotes a quantity which equals $1$ when the event $E$ is true, and $0$ otherwise.

By $\|A\|_{\iy}$ or simply $\|A\|$ we denote the operator norm of a matrix $A$. 


\subsection{Semicircular and Mar\v{c}enko--Pastur  distributions}

Let $m \in \R$ and $\sigma>0$. The {\itshape semicircular distribution with mean $m$ and variance $\sigma^2$} is the probability
distribution $\mu_{SC(m,\sigma^2)}$ with support $[m-2\sigma,m+2\sigma]$ and density
\[ \frac{d\mu_{SC(m,\sigma^2)}}{dx} = \frac{1}{2\pi\sigma^2} \sqrt{4\sigma^2-(x-m)^2}.\]
It is well-known (\cite{agz}, page 7) that if $X$ is a random variable with $SC(0,1)$ distribution, the moments of $X$ are related to the Catalan numbers $C_k=\frac{1}{k+1} \binom{2k}{k}$,
\[ \E X^{2k} = C_k, \ \ \ \ \ \E X^{2k+1} = 0. \]

We now introduce the Mar\v{c}enko--Pastur distributions. First, for $0 < \alpha \leq 1$,
 let $f_\alpha$ be the probability density defined on $[b_-,b_+]$ (where $b_\pm = (1 \pm \sqrt{\alpha})^2)$ by
\[ f_{\alpha}(x) = \frac{\sqrt{(x-b_-)(b_+-x)} }{2\pi x\alpha} .\]
The {\itshape Mar\v{c}enko--Pastur distribution with parameter $\alpha$}, $\mu_{MP(\alpha)}$,
 is the following probability distribution 
\begin{itemize}
 \item If $\alpha \geq 1$, then $\mu_{MP(\alpha)}$ is the probability distribution with density $f_{1/\alpha}$.
 \item If $0 < \alpha \leq 1$, then $d\mu_{MP(\alpha)}(x) = (1-\alpha)\delta_0 + \alpha df_{\alpha}(x)$, where $\delta_0$ denotes a Dirac mass at $0$.
\end{itemize}

In particular, note the following fact: if $X$ has a semicircle $SC(0,1)$ distribution, then $X^2$ has a Mar\v{c}enko--Pastur $MP(1)$ distribution.


\subsection{Asymptotic spectrum of Wishart matrices: Mar\v{c}enko--Pastur distributions} 

\label{subsection:MP}

Define a {\itshape $(n,p)$-Wishart matrix} as a random $n\times n$ matrix $W$ obtained by setting $W=\frac{1}{p}GG^\dagger$, where $G$ is
a $n \times p$ matrix with independent (real or complex\footnote{A complex-valued random variable $\xi$ has a complex $N(0,1)$ distribution if its real and imaginary parts are independent random variables with real $N(0,\frac 12)$ distribution. In particular, $\E |\xi|^2 = 1$.}) $N(0,1)$ entries. The real case and complex case are completely similar. Our
results are valid for both, although only the complex case is relevant to Quantum Information Theory.

Let $A$ be a $n \times n$ Hermitian matrix, and denote $\lambda_1,\dots,\lambda_n$ the eigenvalues of $A$. 
The {\itshape empirical eigenvalue distribution} of $A$, denoted $N_A$, is the probability measure on Borel subsets 
of $\R$ defined as
\[ N_A = \frac{1}{n} \sum_{i=1}^n \delta_{\lambda_i} .\]
In other words, $N_A(B)$ is the proportion of eigenvalues that belong to the Borel set $B$. For large sizes, the empirical eigenvalue distribution of a Wishart matrix approaches a {\itshape Mar\v cenko--Pastur distribution}.

\begin{thm}[Mar\v{c}enko--Pastur, \cite{mp}]
Fix $\alpha >0$. For every $n$, let $W_n$ be a $(n,\lfloor \alpha n \rfloor)$-Wishart matrix. Then the empirical eigenvalue 
distribution 
of $W_n$ approaches a Mar\v{c}enko--Pastur distribution $MP(\alpha)$ in the following sense. For every interval $I \subset \R$ and
any $\e>0$,
\[ \lim_{n \to \iy} \P( |N_{W_n}(I) - \mu_{MP(\alpha)}(I)|>\e)=0.\]
\end{thm}

\subsection{Partial transposition} \label{section-defPT}

We now assume that $n=d^2$. One can think of any $n \times n$ matrix $A$ as a block matrix, consisting of $d \times d$ blocks, each
block being a $d \times d$ matrix. The entries of the matrix are then conveniently described using $4$ indices ranging from $1$ to $d$
\[ A = (A_{i,j}^{k,l})_{i,j,k,l} .\]
Here $i$ denotes the block row index, $j$ the block column index, $k$ the row index inside the block $(i,j)$ 
and $l$ the column index inside the block $(i,j)$. 
We can then apply to each block of $A$ the transposition operation. The resulting matrix is denoted $A^\Gamma$ and called the {\itshape partial transposition}\footnote{An explanation for the notation is that $\Gamma$ is ``half'' of the letter T which denotes the usual transposition.}  of $A$. Using indices, we may write
\begin{equation} \label{def-partialtranspose} (A^\Gamma)_{i,j}^{k,l} = A_{i,j}^{l,k} .\end{equation}

Such a block matrix $A$ can be naturally seen as an operator on $\C^d \otimes \C^d$. Indeed, a natural basis in this space is the double-indexed family $(e_i \otimes e_k)_{1 \leq i,k \leq d}$, where $(e_i)$ is the canonical basis of 
$\C^d$. The action of $A$ on this basis is described as
\[ A(e_i \otimes e_k) = \sum_{j,l=1}^d A_{i,j}^{k,l} e_j \otimes e_l .\]
We may identify canonically $\mM(\C^d \otimes \C^d)$ with $\mM(\C^d) \otimes \mM(\C^d)$. Via this identification, the
matrix $A^\Gamma$ coincides with $(\Id \otimes T)(A)$, where $T: \M(\C^d) \to \M(\C^d)$ is the usual transposition map.
The map $T$ is the simplest example of a map which is positive but not completely positive: $A \geq 0$ does not imply $A^\Gamma \geq 0$.

\subsection{Asymptotic spectrum of partially transposed Wishart matrices: non-centered semicircular distribution}

Motivated by Quantum Information Theory, we investigate the following question: what does the spectrum of $A^\Gamma$ look like ? 
As we will see, the partial transposition dramatically changes the spectrum: the empirical eigenvalue distribution of $A^\Gamma$ is
no longer close to a Mar\v{c}enko--Pastur distribution, but to a shifted semicircular distribution  ! This is our main
theorem.

\begin{N-thm} \label{main-theorem}
Fix $\alpha>0$. For every $d$, let $W_d$ be a $(d^2,\lfloor \alpha d^2 \rfloor)$-Wishart matrix, 
and let $Y_d=W_d^\Gamma$ be the partial transposition of $W_d$. Then the empirical eigenvalue distribution of $Y_d$ approaches the semicircular distribution $\mu_{SC(1,1/\alpha)}$ in the following sense. For every interval
$I \subset \R$ and any $\e >0$,
\[ \lim_{d \to \iy} \P \left( \left| N_{Y_d}(I) - \mu_{SC(1,1/\alpha)}(I) \right| > \e \right) = 0 .\]
Recall that $N_{Y_d}(I)$ is the proportion of eigenvalues of the matrix $Y_d$ that belong to the interval $I$.
\end{N-thm}

Note that the trace and the Hilbert--Schmidt norm are obviously invariant under partial transpose. The distributions  $MP(\alpha)$ and $SC(1,1/\alpha)$ (corresponding to eigenvalue distribution before and after applying partial transpose) indeed share the same first and second moments. 

The support of the limiting spectral distribution $SC(1,1/\alpha)$ is the interval $[1-\frac{2}{\sqrt{\alpha}},1+\frac{2}{\sqrt{\alpha}}]$. Denote by $\lambda_{\min}(A)$ (resp. $\lambda_{\max}(A)$) the smallest (resp. largest) eigenvalue of a matrix $A$. A natural (and harder) question is whether the extreme eigenvalues of $Y_d$ converge towards $1\pm \frac{2}{\sqrt{\alpha}}$. We show that this is indeed the case:

\begin{N-thm} \label{thm:extreme-egv}
Fix $\alpha>0$.  For every $d$, let $W_d$ be a $(d^2,\lfloor \alpha d^2 \rfloor)$-Wishart matrix, 
and let $Y_d=W_d^\Gamma$ be the partial transposition of $W_d$. Then, for every $\e > 0$,
\[ \lim_{d \to \iy} \P \left( \left| \lambda_{\max}(Y_d) - (1+{2}/{\sqrt{\alpha}}) \right| > \e \right) = 0 ,\]
\[ \lim_{d \to \iy} \P \left( \left| \lambda_{\min}(Y_d) - (1-{2}/{\sqrt{\alpha}}) \right| > \e \right) = 0 .\]
\end{N-thm}

\subsection{Almost sure convergence}

In Random Matrix Theory, it is customary to work with the stronger notion of almost sure convergence. This requires to define all the objects on a single probability space. Such a construction is not natural from a Quantum Information Theory point of view, which usually ``avoids infinity'' and prefers to work in a fixed (but large) dimension. 

However, from a mathematical point of view, it is interesting  to note that the results presented here also hold for almost sure convergence. One needs to check that the proof gives enough concentration in order to use the Borel--Cantelli lemma. A key point is the $O(1/d^2)$ estimate for the variance from Proposition \ref{prop:main-variance}.

\section{Non-crossing partitions and combinatorics of Wishart matrices} \label{sec:wishart}

\subsection{Non-crossing partitions}

Let $S$ be a finite set with a total order $<$. Usually, $S$ equals $[k]$ (the set $\{1,\dots,k\}$) for some positive 
integer $k$, and additions in $[k]$ are understood modulo $k$. It is useful to represent elements of $S$ as points on a circle. We introduce the concept of non-crossing partitions and refer to \cite{ns} for more information and pictures.
\begin{itemize}
 \item  A {\itshape partition} $\pi$ of $S$ is a family $\{V_1,\dots,V_p\}$ of disjoint nonempty subsets of $S$, whose union is $S$. The sets $V_i$ are called the {\itshape blocks} of $\pi$. The number of blocks of $\pi$ is denoted $|\pi|$. We denote $\sim_\pi$ the equivalence relation on $S$ induced by $\pi$: $i \sim_\pi j$ means that $i$ and $j$ belong to the same block.
 \item A partition $\pi$ of $S$ is said to be {\itshape non-crossing} if there does not exist elements $i<j<k<l$ in $S$
such that  $i \sim_\pi k, j \sim_\pi l$ and $i \not\sim_\pi j$. We denote by $NC(S)$ the set of non-crossing partitions of $S$, and $NC(k)=NC([k])$.
 \item A {\itshape chording} (or a {\itshape non-crossing pair partition}) of $S$ is a non-crossing partition of 
$S$ in which each block contains exactly two elements. Chordings exist only when the cardinal of $S$ is even. We denote
by $NC_2(S)$ the set of chordings of $S$, and $NC_2(k)=NC_2([k])$. 
\end{itemize}

Counting non-crossing partitions is a well-known combinatorial problem involving Catalan numbers (see \cite{ns}, Lemma 8.9 and Proposition 9.4).

\begin{N-lem} \label{chordings-catalan}
Let $k \in \N^*$. The number of elements in $NC(k)$ and the number of elements in $NC_2(2k)$ are both equal to the 
$k$th Catalan number $C_k=\frac{1}{k+1} \binom{2k}{k}$.
\end{N-lem}

Let us also introduce the {\itshape Kreweras complementation} as the map $K : NC(k) \mapsto NC(k)$ defined as follows. For $ \pi \in NC(\{1^-,\dots,k^-\}) \simeq NC(k)$, $K(\pi)$ is defined as the coarsest partition $\sigma \in NC(\{1^+,\dots,k^+\})\simeq NC(k)$ such that $\pi \cup \sigma$ is a non-crossing partition of $\{1^-,1^+,\dots,k^-,k^+\}$, equipped
with the order
\[ 1^- < 1^+ < 2^- < 2^+ < \cdots < k^- < k^+. \]

The map $K$ is bijective. Moreover,  given $\sigma \in NC(\{1^+,\dots,k^+\})\simeq NC(k)$, one can recover $K^{-1}(\sigma)$ as the coarsest partition $\pi \in NC(\{1^-,\dots,k^-\})\simeq NC(k)$ such that $\pi \cup \sigma$ is a non-crossing partition of $\{1^-,1^+,\dots,k^-,k^+\}$. See \cite{ns} for more details.

The following lemma will be used in connection to partial transposition.

\begin{N-lem} \label{lem:kreweras}
Let $\pi \in NC(k)$ a non-crossing partition and $K(\pi)$ its Kreweras complement. Then,
\begin{enumerate}
 \item For every index $i \in [k]$, 
 \[  \textnormal{ The singleton $\{i\}$ is a block in }K(\pi) \iff i \sim_\pi i+1 .\]
 \item For every distinct indices $i,j \in [k]$,
 \[  \textnormal{ The pair $\{i,j\}$ is a block in }K(\pi) \iff i \sim_\pi j+1 \textnormal{ and } i+1 \sim_\pi j 
 \textnormal{ and } i \not\sim_\pi j  .\]
\end{enumerate}
\end{N-lem}

\begin{proof}
This is geometrically obvious.
\end{proof}

\subsection{Combinatorics related to Wishart matrices}

We now remind the reader about the (standard) proof of the
Mar\v{c}enko--Pastur theorem via the moments method. This proof can be found for example in \cite{jonsson, op} or the book \cite{bs}. Not only our proof will mimic this one, but we will actually strongly recycle
most of the combinatorial lemmas. Let $W_n=(W_{ij})$ be a $(n,p)$-Wishart matrix, and $k \in \N$.
The expansion of $\E \frac{1}{n} \tr W_n^k$ reads
\begin{eqnarray} \label{sum-wishart} 
\nonumber \E \frac{1}{n} \tr W_n^k & = & \frac{1}{n} \sum_{\vec a \in [n]^k} \E W_{a_1,a_2} W_{a_2,a_3} \cdots W_{a_k,a_1} \\
 & = & \frac{1}{np^k} \sum_{\vec a \in [n]^k, \vec c \in [p]^k} 
\E G_{a_1,c_1}\overline{G_{a_2,c_1}}G_{a_2,c_2}\overline{G_{a_3,c_2}} \cdots G_{a_k,c_k}\overline{G_{a_1,c_k}}.
\end{eqnarray}

The next task is to analyze which couples $(\vec a,\vec c)$ give dominant contributions to the sum 
\eqref{sum-wishart} when $n \to \iy$ and $p=\lfloor \alpha n\rfloor$. One argues as follows. First, if one couple $(a_i,c_i)$ or $(a_{i+1},c_i)$ appears a odd number of times in the product, then the contribution is exactly zero (because entries of $G$ are independent and symmetric). This motivates the following definition:

\begin{defn} 
A couple $(\vec a,\vec c) \in \N^k \times \N^k$ satisfies the {\itshape Wishart matching condition} if every couple in the following list of $2k$ elements appears an even number of times:
\begin{equation} \label{list-wishart}  (a_1,c_1), (a_2,c_1), (a_2,c_2), (a_3,c_2), \dots, (a_k,c_k), (a_1,c_k). \end{equation}
\end{defn}

Let $(\vec a,\vec c) \in \N^k \times \N^k$. We define $d_W(\vec a,\vec c)$ as the number of distinct couples appearing in the list \eqref{list-wishart}, and set $\ell_W(\vec a,\vec c) = \# \vec a + \# \vec c$. We also denote $n_2(\vec a,\vec c)$ the number of indices $i$ such that the $i$th element appears exactly twice in the list \eqref{list-wishart}, and 
$n_+(\vec a,\vec c)$ the number of indices $i$ such that the $i$th element appears at least $4$ times. Note that $n_2(\vec a,\vec c)+n_+(\vec a,\vec c)=2k$. These parameters satisfy some inequalities:

\begin{N-lem} \label{lem-wmc}
Let $(\vec a,\vec c) \in \N^k \times \N^k$ satisfy the Wishart matching condition. Then 
\[ \ell_W(\vec a,\vec c) \leq d_W(\vec a, \vec c)+1 \leq k+1. \] 
Moreover, $n_+(\vec a,\vec c) \leq 4(k+1-\ell_W(\vec a,\vec c))$. 
\end{N-lem}

\begin{proof}
Read the list \eqref{list-wishart} from left to right, and count how many new indices you read. The first couple $(a_1,c_1)$ brings two new indices, and each subsequent couple that did not appear previously in the list (there are $d_W(\vec a,\vec c)-1$ such couples) may bring at most one new index (since it shares a common index with the couple just before). This shows that $\ell_W(\vec a,\vec c) \leq d_W(\vec a, \vec c)+1$.

The inequality $d_W(\vec a,\vec c) \leq k$ is easy: if every couple in the list \eqref{list-wishart} appears at least twice, then this list contains at most $k$ different couples.

For the last claim, note that 
\[ d_W(\vec a,\vec c) \leq \frac{1}{2} n_2(\vec a,\vec c) + \frac{1}{4} n_+(\vec a,\vec c) = k -\frac{1}{4}n_+(\vec a,\vec c),\]
with equality iff no element in the list \eqref{list-wishart} appears $6$ times or more.
\end{proof}

Now, the couples $(\vec a,\vec c)$ satisfying  $\ell_W(\vec a,\vec c) < k+1$ are easily shown to have a contribution to the sum \eqref{sum-wishart} which is asymptotically zero. Let us say that $(\vec a,\vec c)$ is {\itshape Wishart-admissible} if it satisfies the matching condition, together with
the equality $\ell_W(\vec a,\vec c) = k+1$. 

If $\vec a \in \N^k$, the {\itshape partition induced by $\vec a$}, denoted $\pi(\vec a)$,
 is the partition
of $[k]$ defined as follows: $i$ and $j$ belong to the same block if and only if $a_i=a_j$.
We say that $\vec a, \vec b \in \N^k$ are {\itshape equivalent}
 $(\vec a \sim \vec b)$ if $\pi(\vec a)=\pi(\vec b)$. 
Similarly, a couple $(\vec a,\vec c)$ is equivalent to a couple $(\vec a\,',\vec c\,')$ if $\vec a \sim \vec a\,'$ and
 $\vec c \sim \vec c\,'$. The next proposition (see \cite{jonsson} or \cite{op} for details) characterizes the combinatorial structure of (equivalence classes of) Wishart-admissible couples.

\begin{N-prop} \label{prop-wishart-admissible}
For every integer $k$, 
\begin{enumerate}
 \item[(a)] 
If $(\vec a,\vec c) \in \N^k \times \N^k$ is Wishart-admissible, then
\begin{enumerate}
 \item[(i)] Each couple in the list \eqref{list-wishart} appears exactly twice. 
One occurrence is of the form $(a_i,c_i)$ while the other occurrence is of the form $(a_{i+1},c_i)$. Moreover, the pair-partition of $[2k]$ induced by the list \eqref{list-wishart} is non-crossing. 
 \item[(ii)] The partitions $\pi(\vec a)$ and $\pi(\vec c)$ are non-crossing, and Kreweras-complementary: $\pi(\vec c)=K(\pi(\vec a))$. In particular, $\vec a$ is determined by $\vec c$ up to equivalence.
\end{enumerate}
\item[(b)] The mapping $(\vec a,\vec c) \mapsto \pi(\vec c)$ induces a bijection between the set of equivalence classes 
of Wishart-admissible  couples in $\N^k \times \N^k$ and the set $NC(k)$.
\end{enumerate}
\end{N-prop}

\begin{eg}
Let us give an example of a Wishart-admissible couple for $k=4$. Let $\vec a = (1,2,2,3)$ and $\vec c=(7,3,7,7)$. Then $\ell_W(\vec a,\vec c)=5$. The list \eqref{list-wishart} reads as
\[ (1,7) ; (2,7) ; (2,3) ; (2,3) ; (2,7) ; (3,7) ; (3,7) ; (1,7). \]
Indeed, each couple appears exactly twice. The partition
induced by this list is
\[ \{ \{1,8\}, \{2,5\}, \{3,4\}, \{6,7\} \}\]
while the partitions induced by $\vec a$ and $\vec c$ are
\[ \pi(\vec a) = \{ \{1\}, \{2,3\} , \{4\} \},\]
\[ \pi(\vec c)= K(\pi(\vec a)) = \{ \{ 1,3,4\}, \{2\} \}. \]
\end{eg}

From Proposition \ref{prop-wishart-admissible}, it is easy to check (if $p \sim \alpha n$) that $\lim_{n \to \iy} \E \frac{1}{n} \tr W_n^k$ coincides with the $k$th moment of the Mar\v cenko--Pastur distribution with parameter $1/\alpha$. To obtain more information than convergence in expectation, one usually needs also a control of the variance of $\frac{1}{n} \tr W_n^k$. The next lemma is then relevant. Actually, the stronger conclusion $\ell_W(\vec a,\vec c) + \ell_W(\vec{a}\,',\vec{c}\,') \leq 2k$ holds, but we do not need this sophistication here.

\begin{N-lem} \label{lemma-wishart-variance}
Let $(\vec a,\vec c)$ and  $(\vec {a}\,'\vec{c}\,')$ be two couples in $\N^k \times \N^k$ satisfying the following
conditions
\begin{enumerate}
 \item[(i)] Each couple in the following list of $4k$ elements appears at least twice:
\begin{equation} \label{list-wishart-variance} 
 (a_1,c_1), (a_2,c_1), \dots, (a_k,c_k), (a_1,c_k) \, ;  \,
(a'_1,c'_1), (a'_2,c'_1),  \dots, (a'_k,c'_k), (a'_1,c'_k)
. \end{equation}
 \item[(ii)] At least some couple appears both in the left half and in the right half of
the list \eqref{list-wishart-variance}.
\end{enumerate}
Then $\ell_W(\vec a,\vec c) + \ell_W(\vec{a}\,',\vec{c}\,') \leq 2k+1$.
\end{N-lem}

\begin{proof}
As before, we read the list \eqref{list-wishart-variance} and keep track of the number of indices. We first read the left half of the list in its natural order. We then read the right half of the list, starting by an element which already appeared in the left half and reading from left
to right---with the convention that $(a'_1,c'_1)$ stands at the right of $(a'_1,c_k)$.

The first element $(a_1,c_1)$ brings two new indices, and each subsequent new couple (there are at most $2k-1$ many, since each couple in the list appears at least twice) brings at most one new index.
\end{proof}

If we want to prove estimates on the extreme eigenvalues of Wishart matrices, we also have to analyze lower-order contributions. We here follow the terminology from \cite{geman}. Let $(\vec a,\vec c) \in \N^k \times \N^k$ satisfy the Wishart matching condition. The elements from the list \eqref{list-wishart} fall into one of the following categories.
\begin{enumerate}
 \item[type 1:] innovations for $\vec a$.
 \item[type 2:] innovations for $\vec c$.
 \item[type 3:] first repetitions of an innovation.
 \item[type 4:] other elements.
\end{enumerate}
The $i$th element in the list \eqref{list-wishart} is an {\itshape innovation} if it contains one index which did  not appear already in the list. When $i=2p$ is even, the $i$th element 
is an {\itshape innovation for $\vec a$} if $a_{p+1} \not\in \{a_j \st j<p\}$. When $i=2p-1$ is odd, the $i$th element is an {\itshape innovation for $\vec c$} if $c_p \not\in \{c_j \st j<p\}$. In particular, the first element of the list \eqref{list-wishart} is always an innovation for $\vec c$.

The $i$th element is the {\itshape first repetition} of an innovation if there is a unique $j<i$ such that the $j$th element from the list \eqref{list-wishart} equals the $i$th element, and moreover this $j$th element is an innovation. 

The following lemma asserts that there are few different couples satisfying the Wishart matching condition which have the same types of elements at the same positions. We refer to \cite{geman} for a proof.

\begin{N-lem} \label{lem:few-type4-elements}
Let $T=(t_1,\dots,t_{2k}) \in \{1,2,3,4\}^{2k}$, and let $U = \card \{ i \in [2k] \st t_i=4 \}$. Say that $(\vec a,\vec c)$ is of type $T$ if, for every $i \in [2k]$, the $i$th element in the list \eqref{list-wishart} has type $t_i$. Then, the number of equivalence classes of couples satisfying the Wishart matching condition which are of type $T$ is bounded by $k^{3U}$.
\end{N-lem}

\subsection{Diagonal elements of Wishart matrices are close to $1$}

We will use the following simple fact in our proof.

\begin{N-lem} \label{lem:diagonal-wishart}
Let $W=(W_{ij})$ be a $(n,p)$-Wishart matrix. Then, for any $\e \in (0,1)$, we have
\[ \P \left( 1-\e \leq \inf_{1 \leq i \leq n} W_{ii} \leq \sup_{1 \leq i \leq n} W_{ii} \leq 1+\e \right) \geq 1-Cn\exp(-cp\e^2) ,\]
where $C,c >0$ are absolute constants.
\end{N-lem}

\begin{proof}
Recall that $W = \frac{1}{p} GG^\dagger$, where $G=(G_{ij})$ is a $n \times p$ matrix with independent $N(0,1)$ entries, so that the diagonal terms of $W_n$ follow a $\chi^2$ distribution
\[ W_{ii} = \frac{1}{p} \sum_{j=1}^p |G_{ij}|^2 .\]
The next fact shows that such distributions enjoy very strong concentration properties.
\begin{N-fact} \label{fact:chi2}
Let $g_1,\dots,g_p$ denote independent (real or complex) $N(0,1)$ random variables, and $X$ be the Euclidean norm of the vector $(g_1,\dots, g_p)$. Then for every $t>0$,
\[ \P \left( | X -\sqrt{p} | > t \right) \leq C' \exp(-c't^2) .\]
\end{N-fact}

Fact \ref{fact:chi2} can be proved by direct calculation or follows from concentration of measure (see e.g.\cite{ledoux}). Indeed, the Euclidean norm is a $1$-Lipschitz function and the expectation of $X$ satisfies the inequalities $\sqrt{p-1} \leq \E X \leq \sqrt{p}$. Lemma \ref{lem:diagonal-wishart} follows from Fact \ref{fact:chi2} via the union bound.
\end{proof}

\section{Proof of Theorem \ref{main-theorem}} \label{sec:overview}

For an integer $d$ and $p=\lfloor \alpha d^2 \rfloor$, let $G_d$ be a $d^2 \times p$ matrix with independent $N(0,1)$ entries, $W_d=\frac{1}{p}G_dG_d^\dagger$ and $Y_d=W_d^\Gamma$.
We denote the entries of $G_d$ as $(G_{i,j}^k)$, where $(i,j) \in [d] \times [d]$ denote the row indices and $k \in [p]$ denotes 
the column index. We label the entries of $W_d$ and $Y_d$ as $(W_{i,i'}^{j,j'})$ and $(Y_{i,i'}^{j,j'})$, where 
$(a,a',b,b') \in [d]^4$ according to the convention described in Section \ref{section-defPT}. 

We have to show that $N_{Y_d}$, the empirical eigenvalue distribution of $Y_d$, approaches a non-centered semicircular distribution
$SC(1,\alpha)$. To handle a more symmetric situation (involving a {\itshape centered} semicircular distribution), we will rather consider $Y_d-\Id$. 
By Lemma \ref{lem:diagonal-wishart}, this matrix is very close to $Z_d=Y_d-\diag(Y_d)$. The latter behaves in a nicer way with respect to 
moments combinatorics. We label the entries of $Z_d$ as $(Z_{i,i'}^{j,j'})_{i,i',j,j' \in [d]}$. We have
\[ Z_{i,i'}^{j,j'} = Y_{i,i'}^{j,j'} {\bf 1}_{(i,j) \neq (i',j')} .\]

The following proposition is central to our work. We defer the proof (the combinatorial part of the moments method) to the next section. 

\begin{N-prop} \label{prop:main-expectation}
For every fixed integer $k$, we have
\[ \lim_{d \to \iy} \E \left( \frac{1}{d^2} \tr (Z_d^k)  \right) =
\begin{cases} \alpha^{-k/2} C_{k/2} & \textnormal{if $k$ is even,} \\ 
                                                   0  & \textnormal{otherwise.} 
                                                   \end{cases} \]
\end{N-prop}

We also show that the variance goes to zero---this is actually simpler.

\begin{N-prop} \label{prop:main-variance}
For every fixed integer $k$, we have
\[ \lim_{d \to \iy} \Var \left( \frac{1}{d^2} \tr (Z_d^k) \right) = 0. \]
\end{N-prop}

The proofs of Proposition \ref{prop:main-expectation} and \ref{prop:main-variance} appear in Sections \ref{sec:expectation} 
and \ref{sec:variance}, respectively.

\begin{proof}[Proof of Theorem \ref{main-theorem} (assuming Propositions \ref{prop:main-expectation} and \ref{prop:main-variance})]
We claim that for any interval $I \subset \R$ and $\e >0$, 
\begin{equation} \label{limit-for-Zd}
 \lim_{d \to \iy} \P \left( \left| N_{Z_d}(I) - \mu_{SC(0,1/\alpha)}(I) \right| > \e \right) = 0 .
\end{equation}

Deriving this from Propositions \ref{prop:main-expectation} and \ref{prop:main-variance} is a completely standard procedure. We only sketch the proof and refer to \cite{agz} (pp 10-11) for more details.
Recall that the Catalan numbers $C_k$ satisfy $C_k \leq 4^k$, and that the support of the $SC(0,1/\alpha)$ distribution is $[-2/\sqrt{\alpha},2/\sqrt{\alpha}]$. We first check that the proportion of eigenvalues outside $J=[-3/\sqrt{\alpha},3/\sqrt{\alpha}]$ is asymptotically zero. For every $\e>0$ and even integer $k$,
\begin{eqnarray*}
\limsup_{d \to \iy} \P\left( N_{Z_d}(J^c) > \e \right) & \leq & \frac{1}{\e} \limsup_{d \to \iy} \E N_{Z_d}(J^c) \\
& \leq & \frac{1}{\e} \limsup_{d \to \iy} \E \int x^k(\sqrt{\alpha}/3)^k dN_{Z_d} \\
& \leq & \frac{1}{\e} (\sqrt{\alpha}/3)^k C_{k/2}\alpha^{-k/2} \\
& \leq & \frac{1}{\e} (2/3)^k ,
\end{eqnarray*}
where the second inequality follows from 
${\bf 1}_{J^c}(x) \leq x^k (\sqrt{\alpha}/3)^k$. Since $k$ is arbitrarily large, we obtain that $\P(N_{Z_d}(J^c)>\e)$
tends to $0$.

Therefore, to prove \eqref{limit-for-Zd}, we may assume $I \subset J$. Using the Weierstrass approximation theorem, we may find a polynomial $Q \geq {\bf 1}_I$ such that $\int Q d\mu_{SC(0,1/\alpha)} \leq \mu_{SC(0,1/\alpha)}(I) + \e/2$. It follows
from Proposition 4.1 that 
\[ \lim_{d \to \iy} \E \int Q dN_{Z_d} = \int Q d\mu_{SC(0,1/\alpha)} .\]
\[ \lim_{d \to \iy} \Var \int Q dN_{Z_d} = 0.\]
For $d$ large enough, $|\E \int Q dN_{Z_d} - \int Q d\mu_{SC(0,1/\alpha)}|<\e/4$. Then
\begin{eqnarray*}
 \P\left(N_{Z_d}(I) \geq \mu_{SC(0,1/\alpha)}(I) +\e\right) &\leq & \P \left( \int Q dN_{Z_d} \geq  \E \int QdN_{Z_d}
 + \e/4 \right) \\
 & \leq & \frac{16}{\e^2} \Var \int Q dN_{Z_d}
 \end{eqnarray*}
and this quantity tends to zero. This is only half of \eqref{limit-for-Zd}. The other half follows by noticing that
\[ \P\left(N_{Z_d}(I) \leq \mu_{SC(0,1/\alpha)}(I) - \e\right) \leq \P \left( N_{Z_d}(J \setminus I) \geq \mu_{SC(0,1/\alpha)}(J \setminus I) +\e/2\right) + \P \left( N_{Z_d}(J^c) \geq \e/2 \right) \]
and applying the previous argument to $J \setminus I$.

We now argue that the empirical eigenvalue distribution is stable under small perturbations. Indeed, for any interval $[a,b]$ and any self-adjoint matrix $\Delta_d$ with operator norm smaller than $\delta$,
\begin{equation} \label{pertubation} N_{Z_d+\Delta_d}([a+\delta,b-\delta]) \leq N_{Z_d}([a,b]) \leq N_{Z_d+\Delta_d}([a-\delta,b+\delta]) .\end{equation}
This is a consequence of the minimax formula for eigenvalues (see e.g. \cite{bhatia}, Chapter III).
We apply \eqref{pertubation} with $\Delta_d=\diag(Y_d)-\Id$. By Lemma \ref{lem:diagonal-wishart}, for every $\e>0$, $\P(\|\Delta_d\|>\e)$ tends to $0$ when $d$ tends to infinity. We easily derive from \eqref{limit-for-Zd} and \eqref{pertubation} that, for any interval $I$,
\[ \lim_{d \to \iy} \P \left( \left| N_{Y_d-\Id}(I) - \mu_{SC(0,1/\alpha)}(I) \right| > \e \right) = 0 . \]
This is clearly equivalent to Theorem \ref{main-theorem}.
\end{proof}

\section{Proof of Proposition \ref{prop:main-expectation}} \label{sec:expectation}

We expand $\E \frac{1}{d^2}\tr(Z_d^k)$ and analyze the underlying combinatorics.
\begin{eqnarray*} \E\frac{1}{d^2} \tr (Z_d^k) &=& \frac{1}{d^2} \sum_{\vec a \in [d]^k, 
\vec b \in [d]^k} \E Z_{a_1,a_2}^{b_1,b_2} \cdot Z_{a_2,a_3}^{b_2,b_3}  \dots  
Z_{a_{k-1},a_k}^{b_{k-1},b_k} \cdot Z_{a_k,a_1}^{b_k,b_1} \\
&=& \frac{1}{d^2} \sum_{\vec a \in [d]^k, 
\vec b \in [d]^k} \left( \prod_{i=1}^k {\bf 1}_{(a_i,b_i) \neq (a_{i+1},b_{i+1})} \right) \E Y_{a_1,a_2}^{b_1,b_2} \cdot Y_{a_2,a_3}^{b_2,b_3} \dots Y_{a_{k-1},a_k}^{b_{k-1},b_k} \cdot Y_{a_k,a_1}^{b_k,b_1} \\
 & = & \frac{1}{d^2} \sum_{\vec a \in [d]^k, \vec b \in [d]^k} M(\vec a,\vec b) \E W_{a_1,a_2}^{b_2,b_1} \cdot W_{a_2,a_3}^{b_3,b_2}  \dots  
W_{a_{k-1},a_k}^{b_{k},b_{k-1}} \cdot W_{a_k,a_1}^{b_1,b_k} \\
& = & \frac{1}{d^2p^k} \sum_{\vec a \in [d]^k, \vec b \in [d]^k, \vec c \in  [p]^k} M(\vec a,\vec b) \E \Pi(\vec a,\vec b,\vec c).
\end{eqnarray*}
where we have defined
\[M(\vec a,\vec b) =  \prod_{i=1}^k {\bf 1}_{(a_i,b_i) \neq (a_{i+1},b_{i+1})}  \]
and
\[ \Pi(\vec a,\vec b, \vec c) =  G_{a_1,b_2}^{c_1} \overline{G_{a_2,b_1}^{c_1}} \cdot G_{a_2,b_3}^{c_2} \overline{G_{a_3,b_2}^{c_2}} \dots G_{a_{k-1},b_k}^{c_{k-1}} 
\overline{G_{a_k,b_{k-1}}^{c_{k-1}}} \cdot G_{a_k,b_1}^{c_k}\overline{G_{a_1,b_k}^{c_k}} .\]

We introduce some definitions in order to restrict ourselves to triples for which both $M(\vec a,\vec b)$ and 
$\E\Pi(\vec a,\vec b,\vec c)$ are nonzero.

\begin{defn}
A couple $(\vec a,\vec b) \in \N^k \times \N^k$ is said to be {\itshape non-repeating} if $M(\vec a,\vec b) = 1$. In other words, $(\vec a,\vec b)$ is
non-repeating if for every $i \in [k]$, either $a_i \neq a_{i+1}$ or $b_i \neq b_{i+1}$.
\end{defn}

Because the entries of $G_d$ are independent, we may factorize $\E \Pi(\vec a,\vec b,\vec c)$ as a product of quantities of the form
$\E (G_{i,j}^k)^q (\overline{G_{i,j}^k})^r$. Such a quantity is zero unless $q=r$, and $\E \Pi(\vec a,\vec b,\vec c)$ is zero whenever 
one of these factors is zero.

\begin{defn} 
A triple $(\vec a,\vec b, \vec c) \in \N^k \times \N^k \times \N^k$ satisfies the {\itshape matching condition} if, in the following list of 
$2k$ triples, each triple appears an even number of times
\begin{equation} \label{liste-PT} (a_1,b_2,c_1), (a_2,b_1,c_1) \, ; 
\, (a_2,b_3,c_2), (a_3,b_2,c_2) \, ; \, \dots \, ; \, (a_k,b_1,c_k), (a_1,b_k,c_k) . \end{equation}
\end{defn}

Therefore, if a triple $(\vec a,\vec b,\vec c)$ does not satisfy the matching condition, then $\E \Pi(\vec a,\vec b,\vec c)=0$ both in the real and in the complex cases. The following easy observation will be used repeatedly.

\begin{N-fact} \label{abc-versus-ac}
Assume that $(\vec a,\vec b,\vec c)$ satisfies the matching condition. Then both $(\vec a,\vec c)$ and $(\vec b,\vec c)$
satisfy the Wishart matching condition.
\end{N-fact}

Recall the definition of equivalence introduced just before Proposition \ref{prop-wishart-admissible}: 
$\vec a \sim \vec a\,'$ means that $\vec a$ and $\vec a\,'$ induce the same partition,
 and $(\vec a,\vec b,\vec c) \sim (\vec a\,', \vec b\,', \vec c\,')$ means
 $\vec a \sim \vec a\,'$, $\vec b \sim \vec b\,'$ and $\vec c \sim \vec c\,'$. Let $C$ be the equivalence class of a triple
 $(\vec a,\vec b,\vec c)$. When $d \to \iy$
\begin{equation}  \label{size-of-eq-class} \card \{C \cap ([d]^k \times [d]^k \times [p]^k)\} \sim d^{\# \vec a} d^{\# \vec b} p^{\# \vec c} \sim 
\alpha^{\# \vec c} d^{\ell(\vec a,\vec b,\vec c)}\end{equation}
where we have defined
\[ \ell(\vec a,\vec b,\vec c) = \# \vec a + \# \vec b + 2\# \vec c . \]
Together with Lemma \ref{lem-wmc}, Fact \ref{abc-versus-ac} implies that whenever $(\vec a,\vec b,\vec c)$ satisfies the matching condition,
\[ \ell(\vec a,\vec b,\vec c) = \ell_W(\vec a,\vec c)+\ell_W(\vec b,\vec c) \leq 2k+2.\]
Let $\mathscr{C}_k$ be the (finite) family of all equivalence classes of triples $(\vec a,\vec b,\vec c) \in \N^k \times \N^k \times \N^k$ which satisfy the matching condition. Since the quantities $M(\vec a,\vec b)$, $\E\Pi(\vec a,\vec b,\vec c)$ and $\ell(\vec a,\vec b,\vec c)$ depend only on the equivalence class $C \in \mathscr{C}_k$ of the triple $(\vec a,\vec b, \vec c)$, we may abusively write $M(C)$, $\E\Pi(C)$ and $\ell(C)$. We also write $\gamma(C)$ to denote $\#\vec c$. Note that these quantities do not 
depend on the dimension $d$. We rearrange the sum according to equivalence classes of triples:
\begin{eqnarray} \label{equation-classes} \lim_{d \to\iy} \frac{1}{d^2} \E \tr Z_d^k & = & \frac{1}{\alpha^k}\sum_{C \in \mathscr{C}_k} M(C) \E\Pi(C)  \lim_{d \to \iy} \frac{1}{d^{2k+2}} \card \{C \cap ( [d]^k \times [d]^k \times [p]^k) \}.
\end{eqnarray}

\begin{N-defn}
Let us say that
 a triple $(\vec a,\vec b,\vec c)$ is {\itshape admissible} if the following three conditions are satisfied
\begin{enumerate}
 \item $(\vec a,\vec b,\vec c)$ satisfies the matching condition,
 \item $(\vec a,\vec b)$ is non-repeating, 
 \item $\ell(\vec a,\vec b,\vec c) = 2k+2$.
\end{enumerate}
Denote by $\mathscr{C}_k^{\textnormal{adm}} \subset \mathscr{C}_k$ the set of equivalence classes of admissible triples.
\end{N-defn}

Equation \eqref{equation-classes} implies that

\begin{equation} \lim_{d \to\iy} \frac{1}{d^2} \E \tr Z_d^k  =  \frac{1}{\alpha^k}\sum_{C \in \mathscr{C}_k^{\textnormal{adm}}} M(C) \E\Pi(C) \alpha^{\gamma(C)}.
\label{last-equation}
\end{equation}

\begin{N-prop} \label{prop:admissible-triples}
If $(\vec a,\vec b,\vec c) \in \N^k \times \N^k \times \N^k$ is admissible, then 
\begin{enumerate}
\item $M(\vec a,\vec b)=1$,
\item $\E \Pi(\vec a,\vec b,\vec c)=1$
\item  $k$ is even
\item $\# \vec c=k/2$. 
\end{enumerate}
Moreover, the number of equivalence classes of admissible triples in $\N^k \times \N^k \times \N^k$ is equal to the 
Catalan number $C_{k/2}$.
\end{N-prop}

Once Proposition \ref{prop:admissible-triples} is proved, Proposition \ref{prop:main-expectation} is immediate from \eqref{last-equation}.

\begin{proof}[Proof of Proposition \ref{prop:admissible-triples}]
The fact that $M(\vec a,\vec b)=1$ is just a reformulation of the non-repeating condition. We now check that $\E \Pi(\vec a,\vec b,\vec c)=1$. Indeed, since $(\vec a,\vec c)$ is Wishart-admissible, every element in the list \eqref{list-wishart} appears exactly twice, once at an odd position and once at an even position. But the same must be true for the list \eqref{liste-PT}, and therefore $\E \Pi(\vec a,\vec b,\vec c)=1$. To check the last two conditions, we rely on the following lemma

\begin{N-lem} \label{lem:bound-on-sizes}
Let $(\vec a,\vec b,\vec c) \in \N^k \times \N^k \times \N^k$ which satisfies the matching condition and such that $(\vec a,\vec b)$ is non-repeating. Then
\begin{enumerate}
 \item No index in $\vec c$ appears only once, and therefore $\# \vec c \leq \lfloor k/2 \rfloor$,
 \item $\# \vec a + \# \vec b \leq 2(\lfloor k/2 \rfloor +1)$.
\end{enumerate}
\end{N-lem}

\begin{proof}
By contraposition, suppose that some index $c_i$ appears only once in $\vec c$, i.e. that $c_j \neq c_i$ for every $j \neq i$. The matching condition imposes the equality
\[ (a_{i+1},b_i,c_i) = (a_i,b_{i+1},c_i) \]
which in turn implies $(a_i,b_i)=(a_{i+1},b_{i+1})$, contradicting the non-repeating property. For the second part of the lemma, we argue differently according to the parity of $k$
\begin{itemize}
 \item[{\bfseries ($k$ odd)}] Define $(\vec x,\vec y) \in \N^k \times \N^k$ as follows
\[ \vec x =(a_1,a_3,\dots,a_{k-2},a_k,a_2,a_4,\dots,a_{k-1}), \ \ \ \ \ 
  \vec y=(b_2,b_4,\dots,b_{k-1},b_1,b_3,\dots,b_{k-2},b_k).
\]
The matching condition implies that $(\vec x,\vec y)$ is Wishart-admissible. Therefore, by Lemma \ref{lem-wmc}, we have $\# \vec x + \# \vec y \leq k+1$. Since
$\vec x$ (resp. $\vec y$) is a permutation of $\vec a$ (resp. $\vec b$), we have 
\[ \# \vec a + \# \vec b \leq k+1 = 2( \lfloor k/2 \rfloor +1). \]
 \item[{\bfseries ($k$ even)}] Define $(\vec x_1,\vec y_1)$ and $(\vec x_2,\vec y_2) \in \N^{k/2} \times \N^{k/2}$ as follows
\[ \vec x_1 = (a_1,a_3,\dots,a_{k-1}), \ \ \ \ \ \vec y_1 = (b_2,b_4,\dots,b_k) ,\]
\[ \vec x_2 = (a_2,a_4,\dots,a_{k}), \ \ \ \ \ \vec y_2  = (b_3,b_5,\dots,b_{k-1},b_1) .\]
Then both $(\vec x_1,\vec y_1)$ and $(\vec x_2,\vec y_2)$ are Wishart-admissible. Therefore, using Lemma \ref{lem-wmc}, we obtain 
\[ \# \vec a +  \# \vec b \leq \# \vec x_1 + \# \vec x_2 + \# \vec y_1 + \# \vec y_2 \leq 2 (k/2+1). \]
\end{itemize}
In both cases we proved $\# \vec a + \# \vec b \leq 2 (\lfloor k/2 \rfloor+1)$.
\end{proof}

We continue the proof of Proposition \ref{prop:admissible-triples}. If $(\vec a,\vec b,\vec c)$ is admissible, 
Lemma \ref{lem:bound-on-sizes} implies that $2k+2 = \ell(\vec a,\vec b,\vec c) \leq 4\lfloor k/2 \rfloor+2$. 
Therefore, $k$ must be even, and necessarily $\# \vec c =k/2$ and each index in $\vec c$ appears exactly twice.

To prove the last statement in Proposition \ref{prop:admissible-triples}, we are going to show that the following map $\Theta$
\begin{eqnarray*} \mathscr{C}_k^{\textnormal{adm}} & \to &NC_2(k) \\
 (\vec a,\vec b,\vec c) & \mapsto & \pi(\vec c) 
\end{eqnarray*}
is bijective. First, the partition induced by $\vec c$ is indeed a chording of $[k]$ (this partition is non-crossing since $(\vec a,\vec c)$ is Wishart-admissible). Because an element of a Wishart-admissible couple is determined (up to equivalence) by the other one, it follows that the map $\Theta$ is injective. 

We now show that this map is onto. Given a partial chording $\pi \in NC_2(k)$, there is a Wishart-admissible couple $(\vec a,\vec c) \in \N^k \times \N^k$ such that $\pi(\vec c)=\pi$. It remains to check that $(\vec a,\vec a,\vec c)$ is admissible.
\begin{itemize}
 \item {\itshape The couple $(\vec a,\vec a)$ is non-repeating}. Otherwise, one would have $a_i=a_{i+1}$ for some index $i \in [k]$. 
 Since $\pi(\vec c)=K(\pi(\vec a))$, this would imply by Lemma \ref{lem:kreweras} that $\{i\}$ is a block in $\pi(\vec
 c)$, which is not possible if $\pi(\vec c)$ is a chording.
 \item {\itshape The triple $(\vec a,\vec a,\vec c)$ satisfies the matching condition}. Since we already know that $(\vec a,\vec c)$
 satisfies the Wishart matching condition, we have to check the following: whenever $(a_i,c_i)=(a_{j+1},c_j)$, we have
 $a_{i+1}=a_j$. Suppose $(a_i,c_i)=(a_{j+1},c_j)$. Since $(\vec a,\vec a)$ is non repeating, we have $i \neq j$. This implies
 that $\{i,j\}$ must be a block in $\pi(\vec c)$ and the result now follows from the second part of Lemma \ref{lem:kreweras}.
\end{itemize}
Therefore, the map $\Theta$ is bijective, and the cardinal of $\mathscr{C}_k^{\textnormal{adm}}$ equals the cardinal of 
$NC_2(k)$, which by Lemma \ref{chordings-catalan} is the Catalan number $C_{k/2}$.
\end{proof}

\section{Proof of Proposition \ref{prop:main-variance}} \label{sec:variance}

Start with a formula from the previous section
\[ \frac{1}{d^2} \tr (Z_d^k) = \frac{1}{d^2p^k} \sum_{\vec a \in [d]^k, \vec b \in [d]^k, \vec c \in  [p]^k} M(\vec a,\vec b)  \Pi(\vec a,\vec b,\vec c). \]
The covariance of two random variables $X,Y$ is defined as $\Cov(X,Y)=\E(XY)-\E X \cdot \E Y$. We have

\begin{equation} \label{expansion-variance} \Var \frac{1}{d^2} \tr Y_d^k = \frac{1}{d^4p^{2k}} \sum_{\vec a,\vec b, \vec c, \vec a\,', \vec b\,', \vec c\,'} M(\vec a,\vec b)M(\vec a\,',\vec b\,') \Cov( \Pi(\vec a,\vec b, \vec c),\Pi(\vec a\,',\vec b\,', \vec c\,') ) , \end{equation}
where the summation is taken over indices $\vec a,\vec b,\vec a\,',\vec b\,'$ in $[d]^k$, and $\vec c,\vec c\,'$ in $[p]^k$. We first identify the vanishing contributions.

\begin{N-lem} \label{lem-variance}
Let $(\vec a,\vec b, \vec c)$ and $(\vec a\,',\vec b\,', \vec c\,')$ be two triples in $\N^k \times \N^k \times \N^k$ such that
\[  \Cov( \Pi(\vec a,\vec b, \vec c),\Pi(\vec a\,',\vec b\,', \vec c\,') ) \neq 0.\]
Then $\ell(\vec a,\vec b,\vec c) + \ell(\vec a\,',\vec b\,', \vec c\,') \leq 4k+2$.
\end{N-lem}

\begin{proof}
The independence of entries of $G_d$ shows that the following two conditions must hold:
\begin{itemize}
 \item Each couple in the following list of $4k$ elements appears at least twice:
\begin{equation} \label{list-variance} 
 (a_1,b_2,c_1), (a_2,b_1,c_1), \dots, (a_k,b_1,c_k), (a_1,b_k,c_k) \, ; \,
(a'_1,b'_2,c'_1), (a'_2,b'_1,c'_1), \dots, (a'_k,b'_1,c'_k), (a'_1,b'_k,c'_k) .
 \end{equation}
 \item At least some couple appears both in the left half and in the right half of
the list \eqref{list-variance}. Otherwise, the random variables $\Pi(\vec a,\vec b, \vec c)$ and $\Pi(\vec a\,',\vec b\,', \vec c\,')$ would be independent, and their covariance would be zero.
\end{itemize}
As is immediately checked, these conditions imply that $\vec a,\vec c,\vec a\,',\vec c\,'$ satisfy the hypotheses of lemma
\ref{lemma-wishart-variance}. Therefore,
\[ \ell_W(\vec a,\vec c) + \ell_W(\vec a\,',\vec c\,') \leq 2k+1 .\]
Similarly, one may apply Lemma \ref{lemma-wishart-variance} to $\vec b,\vec c,\vec b\,',\vec c\,'$ to obtain 
\[ \ell_W(\vec b , \vec c) + \ell_W( \vec b\,' , \vec c\,') \leq 2k+1 .\]
It remains to add both inequalities.
\end{proof}

We now gather the non-zero terms appearing in the sum \eqref{expansion-variance} according to the equivalence class of 
$(\vec a,\vec b,\vec c,\vec a\,',\vec b\,',\vec c\,')$. The cardinality of the equivalence class of $(\vec a,\vec b,\vec c,\vec a\,',\vec b\,',\vec c\,')$ is bounded by
\[ d^{\# \vec a + \# \vec b + \# \vec a\,' + \#\vec b\,'} p^{\# \vec c + \# \vec c\,'}  = 
O \left( d^{\ell(\vec a,\vec b,\vec c) + \ell(\vec a\,',\vec b\,', \vec c\,')}\right) = 
O \left( d^{4k+2}\right).\]
The overall factor $1/d^4p^{2k} = O(1/d^{4k+4})$ in front of the sum \eqref{expansion-variance} shows that  each class has contribution asymptotically zero. Since the number of equivalence classes depends only on $k$,
this proves the lemma.

\section{Convergence of extreme eigenvalues: proof of Theorem \ref{thm:extreme-egv}} \label{sec:extreme}

Let $G_d$ be a $d^2 \times p$ matrix with independent $N(0,1)$ entries, $W_d=\frac{1}{p}G_dG_d^\dagger$, $Y_d=W_d^\Gamma$ 
and $Z_d=Y_d-\diag(Y_d)$. Assume that $p=\lfloor \alpha d^2 \rfloor$. 

Half of Theorem \ref{thm:extreme-egv} can be deduced from Theorem \ref{main-theorem}. Indeed, for every $\e>0$, let $I$ be the 
interval $[1+2/\sqrt{\alpha} - \e, 1+2/\sqrt{\alpha}]$. Since $\mu_{SC(1,1/\alpha)}(I) > 0$, Theorem $1$ implies that, with probability 
tending to $1$, $N_{Y_d}(I) >0$, which means $\lambda_{\max}(Y_d) \geq 1+2/\sqrt{\alpha}-\e$. A similar argument shows that
$\lambda_{\min}(Y_d) \leq 1-2/\sqrt{\alpha}+\e$ with probability tending to $1$.

To prove the other half of Theorem \ref{thm:extreme-egv} (the hard part), we are going to give an upper bound on
$\E \tr(Z_d^k)$ which holds in any fixed dimension (as opposed to asymptotic estimates from the previous sections).

\begin{N-prop} \label{prop:bound-on-moments}
There is a polynomial $Q$ such that, for any integer $k$,
\[ \E \tr(Z_d^k) \leq (2/p)^k (d+Q(k))^{k+2} (\sqrt{p}+Q(k))^k.\]
\end{N-prop}

Assume for the moment that Proposition \ref{prop:bound-on-moments} is true. We claim that it implies that for every $\e>0$,
\[ \lim_{d \to \iy} \P\left( \| Y_d-\Id \| \geq 2/\sqrt{\alpha} + \e \right) = 0 ,\]
from which Theorem \ref{main-theorem} follows. Indeed, choose $k=k(d)$ an even integer such that $Q(k)=o(d)$ and $\log d=o(k)$.
Then, when $d\to \iy$, Proposition \ref{prop:bound-on-moments} implies
\[ \E \|Z_d\|^k \leq \E \tr(Z_d^k) \leq \left( \frac{2d}{\sqrt{p}} +o(1) \right)^k = \left( \frac{2}{\sqrt{\alpha}} +o(1) \right)^k .\]
Therefore, it follows from Markov's inequality that for every $\e > 0$,
\[ \P \left( \|Z_d\| \geq {2}/{\sqrt{\alpha}} +\e \right) \leq \left( \frac{2}{\sqrt{\alpha}} + o(1)\right)^k
\left( \frac{2}{\sqrt{\alpha}} + \e \right)^{-k} \longrightarrow 0. \]
On the other hand, by Lemma \ref{lem:diagonal-wishart},
\[ \P \left( \|\diag(Y_d)-\Id\| \geq \e \right) \leq d^2\exp(-cp\e^2) \longrightarrow 0. \]
This completes the proof of Theorem \ref{thm:extreme-egv} since
\[ \P \left( \|Y_d-\Id\| \geq 2/\sqrt{\alpha} + \e \right) \leq \P \left( \|\diag(Y_d)-\Id\| \geq \e/2\right )+
\P \left( \|Z_d\| \geq 2/\sqrt{\alpha} + \e/2\right ). \]

\begin{proof}[Proof of Proposition \ref{prop:bound-on-moments}]
Recall the computation from Section \ref{sec:expectation}
\begin{equation} \label{eq:sum} \tr (Z_d^k) = \frac{1}{p^k} \sum_{\vec a \in [d]^k, \vec b \in [d]^k, \vec c \in  [p]^k} M(\vec a,\vec b)  \Pi(\vec a,\vec b,\vec c). \end{equation}
We first give an upper bound on $\E \Pi(\vec a,\vec b,\vec c)$.

\begin{N-lem}
Let $(\vec a,\vec b,\vec c) \in \N^k \times \N^k \times \N^k$ satisfy the matching condition, and denote 
\[ \Delta = 2k+2-\ell(\vec a,\vec b,\vec c) .\]
Note that $\Delta \geq 0$. Then
\begin{enumerate}
 \item The number $N$ of indices $i \in [2k]$ such that the $i$th term in the list \eqref{liste-PT} appears $4$ times or more is bounded by $2\Delta$,
 \item We have $\E \Pi(\vec a,\vec b,\vec c) \leq (C_0k)^{\Delta}$, where $C_0$ is an absolute constant.
\end{enumerate}
\end{N-lem}

\begin{proof}
At least one of the numbers $k+1-\ell_W(\vec a,\vec c)$ and $k+1 -\ell_W(\vec b,\vec c)$ must be smaller than $\Delta/2$, since their sum equals $\Delta$. Without loss of generality, we may assume that $k+1-\ell_W(\vec a,\vec c) \leq \Delta/2$. Then, Lemma \ref{lem-wmc} implies that $n_+(\vec a,\vec c) \leq 2\Delta$. Since $N \leq n_+(\vec a,\vec c)$, the first part of the lemma follows.

For the second part, we use independence to write $\E \Pi(\vec a,\vec b,\vec c)$ as a product of quantities of the form 
$\E (G_{i,j}^k)^{q_1} (\overline{G_{i,j}^k})^{q_2} \leq \E |G_{i,j}^k|^{q_1+q_2}$. If $G$ is a $N(0,1)$ random variable, then $\E |G|^{2n}$ equals $1\cdot 3 \cdot 5 \cdots (2n-1)$ in the real case 
and $n!$ in the complex case. In both cases, for some constant $C_0$,
 \begin{equation} \label{moments-N01} \E |G|^q \begin{cases}
  =1 & \textnormal{ if } q=2, \\
  \leq (C_0\sqrt{q})^q & \textnormal{ if } q >2.
 \end{cases}
 \end{equation}
 Bounding each individual 
 factor according to \eqref{moments-N01} and using $q \leq 2k$ leads to
 \[ \E\Pi(\vec a,\vec b, \vec c) \leq (C_0\sqrt{2k})^{N} \]
and the second part of the lemma follows.
\end{proof}

The number of triples in $[d]^k \times [d]^k \times [p]^k$ equivalent to a given triple $(\vec a,\vec b,\vec c)$ is equal to
\[ d(d-1)\cdots(d-\# \vec a+1) \cdot d(d-1) \cdots (d-\# \vec b+1) \cdot p(p-1) \cdots (p-\# \vec c+1) \leq d^{\# \vec a+
\# \vec b} p^{\# \vec c}.
\]
Therefore, it is convenient to rearrange the sum \eqref{eq:sum} according to the values of $\# \vec a + \# \vec b$ and $\# \vec c$. We denote by $m_{\ell_1,\ell_2}$ the number of equivalence classes of triples $(\vec a,\vec b,\vec c) \in  \N^k \times \N^k \times \N^k$ which satisfy the matching condition, with $(\vec a,\vec b)$ non-repeating, $ \# \vec a + \# \vec b = \ell_1$ and $\# \vec c=\ell_2$. It follows from the analysis above that
\begin{equation} \label{one-more-equation} \E \tr (Y^k) \leq \frac{1}{p^k} \sum_{\ell_1,\ell_2} d^{\ell_1} p^{\ell_2} m_{\ell_1,\ell_2} (C_0k)^{2k+2-\ell_1-2\ell_2}. \end{equation}
By Lemma \ref{lem:bound-on-sizes}, $m_{\ell_1,\ell_2}= 0$ if either $\ell_1>k+2$ or $\ell_2>k/2$. It remains to give a bound on the number $m_{\ell_1,\ell_2}$. This is the content of the following proposition (we postpone the proof to the end of the section).

\begin{N-prop} \label{prop:combinatorics}
There is a polynomial $P$ such that the following holds. Denote by $N_{\Delta}$ the number of equivalence classes
of triples $(\vec a,\vec b,\vec c) \in \N^k \times \N^k \times \N^k$ which satisfy the matching condition, with $(\vec a,\vec b)$ non-repeating and $\ell(\vec a,\vec b,\vec c) =2k+2-\Delta$. We have the bound
\begin{equation} \label{eq:sharp-bound} N_{\Delta} \leq 2^k P(k)^\Delta .\end{equation}
\end{N-prop}

\begin{remark}
The bound given in \eqref{eq:sharp-bound} is quite sharp. Indeed, for $\Delta=0$, it gives $N_0 \leq 2^k$. But $N_0$ is exactly the number of equivalence classes of admissible triples considered in Section \ref{sec:expectation}, where this
number was shown to equal the Catalan number $C_{k/2}$, only slightly smaller that $2^k$.
\end{remark}

We continue the proof of Proposition \ref{prop:bound-on-moments}. We have 
\[ m_{\ell_1,\ell_2} \leq N_{2k+2-\ell_1-2\ell_2} \leq 2^k P(k)^{2k+2-\ell_1-2\ell_2} .\]
Plugging this into \eqref{one-more-equation} and denoting $Q$ the polynomial $Q(k)=C_0kP(k)$, 
\begin{eqnarray*}
\E \tr (Z_d^k) & \leq & \frac{2^k}{p^k} \sum_{\ell_1=2}^{k+2} \sum_{\ell_2=1}^{k/2} d^{\ell_1} p^{\ell_2} Q(k)^{2k+2-\ell_1-2\ell_2} \\ 
 & = & \frac{2^k}{p^k} \left( \sum_{\ell_1=2}^{k+2} d^{\ell_1} Q(k)^{k+2-\ell_1} \right)
 \left( \sum_{\ell_2=1}^{k/2} (\sqrt{p})^{2\ell_2} Q(k)^{k-2\ell_2} \right) \\
 & \leq & (2/p)^k (d+Q(k))^{k+2} (\sqrt{p}+Q(k))^k.
\end{eqnarray*}
This completes the proof of Proposition \ref{prop:bound-on-moments}.
\end{proof}

\begin{proof}[Proof of Proposition \ref{prop:combinatorics}]
 For $(\vec a,\vec b,\vec c) \in \N^k \times \N^k \times \N^k$, let $I=I(\vec a,\vec b,\vec c) \subset [k-1]$ be the subset of indices $i$ such that the following condition holds
 \begin{enumerate}
  \item $a_{i+1} \not\in \{ a_j \st j<i+1 \}$ --- one says that $a_{i+1}$ is an innovation,
  \item $b_{i+1} \not\in \{ b_j \st j<i+1 \}$ --- one says that $b_{j+1}$ is an innovation,
  \item $c_{i} \not\in \{ c_j \st j<i \}$ --- one says that $c_j$ is an innovation.
\end{enumerate}

The next lemma shows that the set $I(\vec a,\vec b,\vec c)$ is large when $\Delta$ is small. We postpone the proof.

\begin{N-lem} \label{lem:many-innovations}
If $(\vec a,\vec b,\vec c) \in \N^k \times \N^k \times \N^k$ satisfies the matching condition with $(\vec a,\vec b)$ non-repeating, then
\[ \card I(\vec a,\vec b,\vec c) \geq k/2 -\Delta. \]
where $\Delta = (2k+2) - \ell(\vec a,\vec b,\vec c)$.
\end{N-lem}

Let $A,C$ be subsets of $[k]$. A couple $(\vec a,\vec c)$ satisfying the Wishart matching condition is said to be compatible with $(A,C)$ if
\begin{enumerate}
 \item for every $i \in A$, the index $a_i$ is an innovation, i.e. $a_i \notin \{a_j \st j<i\}$,
 \item for every $i \in C$, the index $c_i$ is an innovation, i.e. $c_i \notin \{c_j \st j<i\}$.
\end{enumerate}

Note that if a Wishart-admissible couple $(\vec a,\vec c)$ is compatible with $(A,C)$, then by arguing as in the proof of Lemma \ref{lem-wmc}, we have
\[ \card A+\card C \leq d_W(\vec a,\vec c) +1 \leq k+1 .\]

Let us state one more lemma, postponing the proof.

\begin{N-lem} \label{lem:label-compatibility}
Let $A,C$ be subsets of $[k]$, and $\delta = k+1 - \card A - \card C$. The number of equivalence classes 
of couples $(\vec a,\vec c) \in \N^k \times \N^k$ which satisfy the Wishart matching condition and are compatible with 
$(A,C)$ is bounded by $(2k)^{9\delta}$.
\end{N-lem}

The number $N_\Delta$ is the number (up to equivalence) of triples $(\vec a,\vec b,\vec c)$  which satisfies the matching condition, with $(\vec a,\vec b)$ non repeating, and $\ell(\vec a,\vec b,\vec c) = (2k+2) -\Delta$. To bound $N_\Delta$, we first choose a set $I \subset [k-1]$ of cardinal larger than $k/2 - \Delta$. The number of possibilities for $I$ is bounded by $2^k$. Now, given $I$, let $I^+$ be the subset of $[k]$ defined as 
\[ j \in I^+ \iff j=1 \textnormal{ or } j-1 \in I.\]

If $(\vec a,\vec b,\vec c)$ satisfies the matching condition with $I(\vec a,\vec b,\vec c)=I$, then it is easily checked that both couples $(\vec a,\vec c)$ and $(\vec b,\vec c)$ are compatible with $(I^+,I)$. We have $\card (I^+)+\card(I) = 2\card(I) +1 \geq k+1-2\Delta$. By Lemma \ref{lem:label-compatibility}, the number of admissible couples compatible with $(I^+,I)$ is bounded by $(2k)^{18\Delta}$. Therefore the number of possible triples $(\vec a,\vec b,\vec c)$ is bounded by $(2k)^{36}$.
This yields the bound
\[ N_\Delta \leq 2^k (2k)^{36\Delta}. \]
This proves Proposition \ref{prop:combinatorics} with $P(k)=(2k)^{36}$.
\end{proof}

\begin{proof}[Proof of Lemma \ref{lem:many-innovations}]
For each index $i \in [k]$, one of the following possibility occurs
\begin{enumerate}
 \item[$P_1(i)$:] The indices $a_{i+1}, b_{i+1}$ and $c_i$ are innovations. Necessarily the triples $(a_i,b_{i+1},c_i)$ and
 $(a_{i+1},b_{i},c_i)$ are innovations\footnote{We say that a triple at $j$th position from the list \eqref{liste-PT} is an innovation  if it does not coincide with a triple at $i$th position for $i<j$.}.
 \item[$P_2(i)$:] The triples $(a_i,b_{i+1},c_i)$ and $(a_{i+1},b_{i},c_i)$ are innovations, but at least one of
 $a_{i+1},b_{i+1}$ and $c_i$ is not an innovation.
 \item[$P_3(i)$:] Only one of the triples $(a_i,b_{i+1},c_i)$ and $(a_{i+1},b_{i},c_i)$ is an innovation.
 \item[$P_4(i)$:] Neither $(a_i,b_{i+1},c_i)$ nor $(a_{i+1},b_{i},c_i)$ is an innovation.
\end{enumerate}
For $j\in \{1,2,3,4 \}$, let $n_j$ be the number of indices $i \in [k]$ such that $P_j(i)$ holds in the above alternative. With this notation, $n_1=\card I(\vec a,\vec b,\vec c)$. The numbers $n_1,n_2,n_3,n_4$ satisfy the following relations
\begin{equation} \label{system-eq1} n_1+n_2+n_3+n_4 = k, \end{equation}
\begin{equation} \label{system-eq2} n_3+2n_4 \geq k, \end{equation}
\begin{equation} \label{system-eq3} 4n_1+3n_2+n_3 \geq 2k-\Delta. \end{equation}
\begin{itemize}
 \item Equation \eqref{system-eq1} is obvious since possibilities $P_1(i), \dots, P_4(i)$ are mutually exclusive.
 \item There must be at least $k$ elements in the list \eqref{liste-PT} which are not innovations, since every element must 
 appear at least twice. But the number of non-innovations in the list \eqref{liste-PT} is equal to $n_3+2n_4$, hence the 
 equation \eqref{system-eq2}. 
 \item For each $i$, let $Z_i$ be the number
 \[Z_i = {\bf 1}_{\{a_{i+1} \textnormal{ is an innovation}\}} +  {\bf 1}_{\{b_{i+1} \textnormal{ is an innovation}\}}
+ 2 \cdot {\bf 1}_{\{c_i \textnormal{ is an innovation}\}} . \]
The value of $Z_i$ depends on which of $P_1(i),P_2(i),P_3(i),P_4(i)$ occurs. 
If $P_1(i)$ occurs, then $Z_i=4$.
 If $P_4(i)$ occurs, then $Z_i=0$. If $P_2(i)$ occurs, then $Z_i \leq 3$. If $P_3(i)$ occurs, then $Z_i \leq 1$.
 This last point deserves some explanation. 
\begin{itemize}
\item If $(a_i,b_{i+1},c_i)$ is not an innovation, then certainly $b_{i+1}$ and $c_i$
 cannot be innovations. 
\item If instead $(a_{i+1},b_i,c_i)$ is not an innovation, then $a_{i+1}$ cannot be an innovation.
 We claim that $c_i$ is also not an innovation. Indeed, if $c_i$ was an innovation,
 then necessarily $(a_{i+1},b_i,c_i)$ would be equal to $(a_{i},b_{i+1},c_i)$ which would contradict the non-repeating property. 
\end{itemize}
This shows that $\sum Z_i \leq 4n_1+3 n_2+n_3$. On the other hand, we have
\[ \sum_{i=1}^k Z_i = \# \vec a-1+\#\vec b-1 + 2 \#\vec c= 2k-\Delta.\]
Therefore, the above discussion implies equation \eqref{system-eq3}.
\end{itemize}
Adding \eqref{system-eq3} and twice \eqref{system-eq2}, we obtain
\[ 4n_1+3n_2+3n_3+4n_4 \geq 4k-\Delta .\]
Together with \eqref{system-eq1}, this implies that $n_2+n_3 \leq \Delta$. Since $n_3 \geq 0$, this in turn implies $3n_2+n_3 \leq 3\Delta$. Combined with \eqref{system-eq3}, we obtain $4n_1 \geq 2k-4\Delta$, hence $n_1 \geq k/2-\Delta$ as claimed.
\end{proof}

\begin{proof}[Proof of Lemma \ref{lem:label-compatibility}]
Given a couple $(\vec a,\vec c) \in \N^k \times \N^k$ satisfying the Wishart matching condition, there is a partition of $[2k]$ as
\begin{equation} \label{partition} [2k] = T_1 \cup T_2 \cup T_3 \cup T_4 \end{equation}
where $T_i$ denotes the set of indices $j$ such that the $j$th element in the list \eqref{list-wishart} is of type $i$ (the four possible types have been defined in Section \ref{sec:background}). If
the couple $(\vec a,\vec c)$ is compatible with $(A,C)$, then necessarily $T_1^* \subset T_1$ and $T_2^* \subset T_2$,
where
\[ T_1^* = \{2(i-1) \st i \in A, i \neq 1 \}, \]
\[ T_2^* = \{2i-1 \st i \in C \}. \]

We claim that the number of partitions \eqref{partition} satisfying these constraints is bounded by $(2k)^{3\delta}$. Indeed,
we first have to enlarge $T_1^*$ into $T_1$ and $T_2^*$ into $T_2$. Since $\card(T_1^* \cup T_2^*) = k-\delta$ and $\card(T_1 \cup T_2) \leq k$, the number of possible ways to perform these enlargements in at most $(2k)^\delta$.

Since $\card(T_3)=\card(T_1)+\card(T_2)$, we have $\card (T_4) \leq 2\delta$. Therefore the number of possible choices for
$T_4$ is bounded by $(2k)^{2\delta}$. Once $T_1,T_2$ and $T_4$ are chosen, the set $T_3$ consists of the remaining indices. Hence the claim on the number of possible partitions.

Now, by Lemma \ref{lem:few-type4-elements}, the number of equivalence classes of couples satisfying the Wishart matching
condition with a given partition \eqref{partition} is bounded by
\[ (2k)^{3 \card T_4} \leq (2k)^{6\delta} .\]

Finally, the total number of equivalence classes satisfying the Wishart matching condition and compatible with $(A,C)$ is 
bounded by $(2k)^{9\delta}$.
\end{proof}

\section{Relevance to Quantum Information Theory} \label{sec:QIT}

In this section we consider finite-dimensional complex Hilbert spaces. We write $\mM(\C^n)$ for the space
of linear operators (=matrices) on $\C^n$. 

\subsection{PPT states}

A {\itshape state} (=density matrix) $\rho $ on $\C^n$ is a positive operator on $\C^n$ with trace $1$. 
We write $\mD(\C^n)$ for the set of states on $\C^n$. A {\itshape pure state} is a rank one state and is denoted $\rho = \ketbra{x}{x}$, where $x$ is a unit vector in the range of $\rho$. We typically consider the case $\C^n \simeq \C^d \otimes \C^d$.  We have the following canonical identification
\[ \mM(\C^d \otimes \C^d) \simeq \mM(\C^d) \otimes \mM(\C^d) .\]
A state $\rho \in \mD(\C^d \otimes \C^d)$ if called {\itshape separable} if it can be written as a convex combination of product states. A state $\rho$ is called PPT (``positive partial transpose'') if $\rho^\Gamma$ is a positive operator (the partial transposition $\rho^\Gamma = (\Id \otimes T)\rho$ was defined in \eqref{def-partialtranspose}). The partial transposition of a {\itshape separable state} $\rho$ is always positive \cite{peres}; however there exist non-separable (=entangled) PPT states. For many purposes, checking positivity of the partial transpose is the most efficient tool to detect entanglement. We refer to the survey \cite{hhhh-survey} for more information about PPT states and entanglement.

\subsection{Random induced states are normalized Wishart matrices}

 There is a canonical probability measure on the set of pure states on any finite-dimensional Hilbert space $H$, obtained by pushing forward the uniform measure on the unit sphere of $H$ under the map $x \mapsto \ket{x}\bra{x}$. 
We define the measure $\mu_{n,p}$ to be the distribution of $\tr_{\C^p} \ket{x}\bra{x}$, where $x$ is uniformly distributed on the unit sphere of $\C^n \otimes \C^p$. The {\itshape partial trace} $\tr_{\C^p}$ is the linear operation 
\[ \tr_{\C^p}: = \Id_{\mM(\C^n)} \otimes \tr : \mM(\C^n \otimes \C^p) \to \mM(\C^n), \]
 where $\Id_{\mM(\C^n)}$ is the identity operation on $\mM(\C^n)$ and $\tr : \mM(\C^p) \to \C$ is the usual trace. 

The measure $\mu_{n,p}$ is a probability measure on $\mD(\C^n)$, the set of mixed states on $\C^n$. A random state with distribution $\mu_{n,p}$ is called an {\itshape induced state}; the space
$\C^p$ is called the {\itshape ancilla space}. This family of measures has a simple physical motivation: they can be used if our only knowledge about a state is the dimensionality of the environment (see \cite{bz-book}, Section 14.5 and references therein). 

Induced states are closely related to Wishart distributions. Indeed, if $W$ is a $(n,p)$-Wishart random matrix, then $\frac{1}{\tr W} W$ is a random state with distribution $\mu_{n,p}$. Moreover, the random variables $\tr W$ and $\frac{1}{\tr W} W$ are independent (this fact explicitly appears in \cite{nechita}). Therefore, results about Wishart matrices can be easily translated in the language of induced states. The special case $p=n$, when the dimension of the ancilla equals the dimension of the system, deserves to be highlighted thanks to the following Proposition \cite{zs}. 

\begin{N-prop} \label{inducing-HS}
The measure $\mu_{n,n}$ is equal to the normalized Lebesgue measure restricted to the set $\mD(\C^n)$.
\end{N-prop}

Proposition \ref{inducing-HS} follows from a more general fact \cite{zs}:  whenever $p \geq n$, the density of the measure $\mu_{n,p}$ with respect to the Lebesgue measure on $\mD(\C^n)$ is proportional to $\det(\rho)^{p-n}$. 

\subsection{Partial transposition of random induced states}

Our main results admit an immediate translation in the language of random induced states. Here is a version of Theorem \ref{main-theorem} for induced states.

\begin{N-thm} \label{thm:main-states}
Fix $\alpha > 0$. For each $d$, let $\rho_d$ be a random state on $\C^d \otimes \C^d$ chosen according to the measure $\mu_{d^2,\lfloor \alpha d^2 \rfloor}$. Then for every interval $I=[a,b] \subset \R$ and $\e >0$,
\[ \lim_{d \to \iy} \P \left( \left| N_{d^2 \rho_d^\Gamma}(I) - \mu_{SC(1,1/\alpha)}(I) \right| > \e \right) = 0 .\]
Recall that $N_{d^2\rho_d^\Gamma}(I)$ is the proportion of eigenvalues of the matrix $\rho_d^\Gamma$ that belong to the interval $[a/d^2,b/d^2]$.
\end{N-thm}

\begin{proof}
If $W$ is a $(d^2,p)$-Wishart matrix, then $\frac{W}{\tr W}$ has distribution as $\mu_{d^2,p}$. Therefore,
\[ N_{d^2 \rho_d^\Gamma}([a,b])  = N_{\frac{d^2}{\tr W}W^\Gamma}([a,b]) = N_{W^\Gamma}\left( \left[ \frac{\tr W}{d^2} a, \frac{\tr W}{d^2} b \right] \right) .\]
The distribution of $\frac{\tr W}{d^2}$ is proportional to a $\chi^2$ distribution. Using Fact \ref{fact:chi2} to quantify its concentration, we obtain that for any $\eta>0$,
\begin{equation} \label{eq:chi2-deviation} \P \left(\left|\frac{\tr W}{d^2} -1 \right| > \eta \right) \leq C\exp(-cd^2p\eta^2).\end{equation}
When $\left|\frac{\tr W}{d^2} -1 \right| \leq \eta$, we may use the inclusions
\[ [(1+\eta)a,(1-\eta)b] \subset \left[\frac{\tr W}{d^2}a,\frac{\tr W}{d^2}b\right] \subset [(1-\eta)a,(1+\eta)b] \]
to show that Theorem \ref{main-theorem} implies Theorem \ref{thm:main-states}.
\end{proof}

If $d$ is fixed, the induced measures $\mu_{d^2,p}$ concentrate towards the maximally mixed state on $\C^d \otimes \C^d$ when $p$ increases. For small values of $p$, one expects to get typically very entangled states. Therefore one
can consider the critical $p$ for which the property ``being PPT'' becomes typically true. The following theorem shows that a threshold occurs when $p=4d^2$.

\begin{N-thm} \label{thm:threshold-states}
For every $\e>0$, there exist positive constants $c(\e),C(\e)$ such that the following holds. If $\rho$ is a random state on $\C^d \otimes \C^d$ chosen according to the measure $\mu_{d^2,p}$, then
\begin{enumerate}
 \item If $p\leq (4-\e)d^2$, then 
 \[ \P( \rho \textnormal{ is PPT}) \leq C(\e) \exp ( -c(\e)p). \] 
 \item If $p \geq (4+\e)d^2$, then
 \[ \P( \rho \textnormal{ is PPT}) \geq 1-C(\e) \exp ( -c(\e)p). \] 
\end{enumerate}
\end{N-thm}

\begin{proof}
We only show the proof of (1), the proof of (2) being similar. We are going to use a concentration argument from \cite{asy}, where the same question is studied for separability instead of PPT. We start by a lemma that compares the probability that a random state is PPT, for different dimensions.

\begin{N-lem} \label{lem:monotonicity}
Let $d_1,d_2,d'_1,d'_2$ and $p$ be integers, with $d'_1 \leq d_1$ and $d'_2 \leq d_2$. Let $\rho$ be a random state on $\C^{d_1} \otimes \C^{d_2}$ with distribution $\mu_{d_1d_2,p}$, and let $\rho'$ be a random state on $\C^{d'_1} \otimes \C^{d'_2}$ with distribution $\mu_{d'_1d'_2,p}$. Then
\[ \P( \rho \textnormal{ is PPT}) \leq \P(\rho' \textnormal{ is PPT}). \]
\end{N-lem}

\begin{proof}
It is enough to prove the lemma in the special case $d_2=d'_2$ (since both factors play the same role, the full version follows by applying twice this special case). 

We construct a coupling between both distributions as follows. Identify $\C^{d'_1}$ as a subspace of $\C^{d_1}$, and let $Q : \C^{d_1} \to \C^{d'_1}$ be the orthogonal projection. Then, $\C^{d'_1} \otimes \C^{d_2} \subset \C^{d_1} \otimes \C^{d_2}$ is the range of the projection $P = Q \otimes \Id$. Let $W$ be a $(d_1d_2,p)$-Wishart matrix, seen as an operator on $\C^{d_1} \otimes \C^{d_2}$. The random operator $PWP$, when seen
as an operator on $\C^{d'_1} \otimes \C^{d_2}$, has the distribution of a $(d'_1d_2,p)$-Wishart matrix. Therefore, the states 
\[ \rho = \frac{W}{\tr W},\]
\[ \rho' = \frac{P \rho P}{\tr P\rho P} = \frac{PWP}{\tr{PWP}}, \]
have respective distributions $\mu_{d_1d_2,p}$ and $\mu_{d_1'd_2,p}$. To prove the lemma it remains to check that
\[ \rho \text{ is PPT} \Longrightarrow \rho' \text{ is PPT}.\]
This implication holds because $(P\rho P)^\Gamma = P\rho^\Gamma P$.
\end{proof}

Fix $\e>0$. As a consequence of Lemma \ref{lem:monotonicity}, it is enough to prove Theorem \ref{thm:threshold-states}, for every given $p$, when $d$ is minimal such that $p\leq (4-\e)d^2$ (from now one, we assume that $d$ and $s$ are related by  this condition). 

Denote by $\|\cdot\|_{\text{PPT}}$ the gauge associated to the convex body of all PPT states. This gauge is defined as follows, for any state $\rho$ on $\C^d \otimes \C^d$
\begin{eqnarray*} \|\rho\|_{\text{PPT}} & = & \inf \left\{ t \geq 0 \st \frac{\Id}{d^2} + \frac{1}{t} \left(\rho-\frac{\Id}{d^2}\right) \text{ is PPT} \right\} \\
 & = & 1-d^2\lambda_{\min} (\rho^\Gamma).
\end{eqnarray*}

Note in particular that $\rho$ is PPT if and only if $\|\rho\|_{\text{PPT}} \leq 1$. Let $\rho_{d^2,p}$ be a random state with distribution $\mu_{d^2,p}$, and denote by $M_{d^2,p}$ the median of the random variable $\|\rho_{d^2,p}\|_{\text{PPT}}$. By applying Proposition 4.2 from \cite{asy}, we obtain the following inequality: there are absolute constants $c,C$ such that for any $\eta>0$,
\begin{equation}
\label{eq:concentration-inequality}
 \P \left( \left| \|\rho\|_{\text{PPT}} - M_{d^2,p} \right| \geq \eta \right) \leq C \exp (-cp) + C \exp(-cp\eta^2). 
 \end{equation}

Let $W_{d^2,p}$ be a $(d^2,p)$-Wishart matrix. It follows from Theorem \ref{thm:extreme-egv} that $\lambda_{\min}(W_{d^2,p}^\Gamma)$ converges in probability towards $1-2/\sqrt{4-\e}$ when $d,p$ tend to infinity. By \eqref{eq:chi2-deviation}, $\tr W_{d^2,p}/d^2$ converges in probability to $1$. 
Since $W_{d^2,p}/\tr W_{d^2,p}$ has distribution $\mu_{d^2,p}$, it follows that $\|\rho_{d^2,p}\|_{\text{PPT}}$ converges to $\frac{2}{\sqrt{4-\e}}$. In particular,
\[ \lim_{p,d \to \iy} M_{d^2,p} = \frac{2}{\sqrt{4-\e}} > 1.\]

We now choose $\eta$ such that $2/\sqrt{4-\e} > 1+\eta$. For $d,p$ large enough, we have $M_{d^2,p} > 1+\eta$, and
we can apply \eqref{eq:concentration-inequality} to obtain
\[ \P(\rho \textnormal{ is PPT}) = \P( \|\rho\|_{\text{PPT}} \leq 1) \leq C \exp (-cp) + C \exp(-cp\eta^2) .\]

This concludes the proof of Theorem \ref{thm:threshold-states} (small dimensions can be taken into account by adjusting the constants).
\end{proof}

\section{Miscellaneous remarks} \label{sec:misc}

\subsection{Partial transposition of a random pure state}

Another natural question from the point of view of Quantum Information Theory is to study the partial transposition
of random {\itshape pure} states (as opposed to random {\itshape mixed} states considered here). In that direction, one may prove the following result

\begin{N-prop} \label{prop:random-pure-state}
For every $d$, let $\rho_d$ be a random pure state on $\C^d \otimes \C^d$, with uniform distribution. Then, when $d$ tends to infinity, the empirical eigenvalue distribution of $d\rho_d^\Gamma$ approaches a deterministic distribution which can be described as the distribution of the product of two independent $SC(0,1)$ random variables. 
\end{N-prop}

\begin{remark}
The notion of convergence used is the same as in Theorem \ref{thm:main-states}. The limiting distribution appearing in Proposition \ref{prop:random-pure-state} has vanishing odd moments and even moments equal to the square of Catalan numbers. Such a distribution has been studied recently in \cite{BFP}, where a closed formula for the density (involving special functions) is derived.
\end{remark}

\begin{proof}[Proof of Proposition \ref{prop:random-pure-state} (sketch)]
If $\rho=\ketbra{\psi}{\psi}$ is a pure state on $\C^d \otimes \C^d$, the eigenvalues of $\rho^\Gamma$ can be described from the Schmidt coefficients of $\psi$ (Schmidt coefficients for tensors correspond to singular values for matrices, and are therefore governed by the Mar\v{c}enko--Pastur distribution).
Indeed, given a Schmidt decomposition 
\[ \psi = \sum_{i=1}^d \sqrt{\lambda_i} e_i \otimes f_i ,\]
for some orthonormal bases $(e_i),(f_i)$, one checks that
\[ \ketbra{\psi}{\psi}^\Gamma = \sum_{i,j=1}^d \sqrt{\lambda_i \lambda_j} \ketbra{e_i \otimes f_j}{e_j \otimes f_i}. \]
It follows that the eigenvalues of $\ketbra{\psi}{\psi}^\Gamma$ are
\[ \begin{cases} \lambda_i & \textnormal{ for every } 1 \leq i \leq d, \\ \pm \sqrt{\lambda_i\lambda_j} & 
\textnormal{ for every } 1 \leq i < j \leq d.    \end{cases} \]
Eigenvalues of the first category do not contribute to the limit distribution, and the result follows with little effort.
\end{proof}

\subsection{Unbalanced bipartite systems} 

We may apply partial transposition to any decomposition
$\C^{d^2} \simeq \C^{d_1} \otimes \C^{d_2}$, with $d_1d_2=d^2$. 
Provided the ratio $d_1/d_2$ stays away from $0$ and $\iy$, Theorems \ref{main-theorem} and \ref{thm:extreme-egv} remain valid. The point is that the main contributions come from terms in which $\vec a \sim \vec b$, so that $d_1^{\# \vec a} d_2^{\# \vec b}$ depends only on the product $d_1d_2$. 

\subsection{Connexions to free probability}

The same model of partially transposed Wishart matrices has been considered recently by Banica and Nechita \cite{bn} in
a different asymptotic regime (when $d_1$ is fixed and $d_2$ goes to infinity). For that regime the picture is different: they obtain that the limit spectral distribution can be described as the difference of two freely independent random variables with Mar\v{c}enko--Pastur distributions. The shifted semicircle distribution appears then as a limit case. We refer to \cite{bn} for more information.

\subsection{Uniform mixtures of random pure states}

There is another popular model of random states which is very similar to the model of random induced states considered in Section
\ref{sec:QIT}, for which our results are also valid. Let $(\psi_i)_{1 \leq i \leq p}$ be unit vectors in $\C^n$,
chosen independently according to the uniform probability measure on the the unit sphere. Then we consider the random state
\[ \rho = \frac{1}{p} \sum_{i=1}^p \ketbra{\psi_i}{\psi_i} .\]
Denote by $\nu_{n,p}$ the distribution of $\rho$. This model of random states has been considered for example in \cite{znidaric}. When $n,p$ are large, the probability measures $\mu_{n,p}$ and $\nu_{n,p}$ behave similarly. It can be shown that Theorems \ref{thm:main-states} and \ref{thm:threshold-states} remain valid when the probability measures $\mu_{n,p}$ are substituted by the probability measures $\nu_{n,p}$. 

\subsection{Volume of the PPT convex body}

How many states have a positive partial transpose ? This question may be formulated using the Lebesgue measure (or ``volume'') induced by the Hilbert--Schmidt scalar product, or equivalently (cf Proposition \ref{inducing-HS}) by the induced measure over an ancilla of equal dimension. Let $W_d$ be a $(d^2,d^2)$-Wishart random matrix. It was shown in \cite{as} (formulated as a lower bound on the volume of the set of PPT states, and using techniques from high-dimensional convexity) that for some constant $C>0$
\begin{equation} \label{volume-PPT} \P(W_d^\Gamma \geq 0) \geq \exp(-C d^4) .\end{equation}
By Theorem \ref{thm:extreme-egv}, the probability on the left-hand side tends to $0$ when $d$ tends to $+\iy$. How fast it goes to zero is actually a question about large deviations. For standard models of random matrices, very precise results are known about large deviations (see e.g. \cite{agz}, Section 2.6.2), and one may expect the lower bound from \eqref{volume-PPT} to be sharp. 

\begin{conj}
There is an absolute constant $c>0$ such that, whenever $W_d$ is a $(d^2,d^2)$-Wishart matrix,
\[ \P(W_d^\Gamma \geq 0) \leq \exp(-c d^4) .\]
\end{conj}

\noindent This would quantify precisely how (un)common are PPT states  in large dimensions.

\end{document}